\documentclass[10pt]{article}
\usepackage{latexsym}
\usepackage{amsfonts}
\usepackage{amsmath}
\usepackage{amssymb}

\newtheorem{thm}{Theorem}

\newtheorem{cnj}[thm]{Conjecture}

\def\Nt{\hat{M}}

\def\qed{\hfill $\vcenter{\hrule height .3mm
\hbox {\vrule width .3mm height 2.1mm \kern 2mm
\vrule width .3mm height 2.1mm} \hrule height .3mm}$ \bigskip}

\newcommand{\R}{\mathbb{R}}
\newcommand{\Rn}{\R^n}
\def\wrt{with respect to }

\begin{document}
\title{On convexified packing and entropy duality
}
\date{}
\author{ S. Artstein,
V. Milman,
S. Szarek,
N. Tomczak-Jaegermann}
\maketitle
\noindent {\bf 1. Introduction.}
If $K$ and $B$ are subsets of a vector space (or just a group, or even
a homogeneous space), the {\em covering number} of $K$ by $B$, denoted
$N(K, B)$, is the minimal number of translates of $B$ needed to cover
$K$. Similarly, the {\em packing} number $M(K, B)$ is the maximal
number of disjoint translates of $B$ by elements of $K$. The two
concepts are closely related; we have $N(K, B-B) \le M(K, B) \le N(K,
(B-B)/2)$. If $B$ is a ball in a normed space (or in an appropriate
invariant metric) and $K$ a subset of that space (the setting and the
point of view we will usually employ), these notions reduce to
considerations involving the smallest $\epsilon$-nets or the largest
$\epsilon$-separated subsets of $K$.

Besides the immediate geometric framework, packing and covering
numbers appear naturally in numerous subfields of mathematics, ranging
from classical and functional analysis through probability theory and
operator theory to information theory and computer science (where a
code is typically a packing, while covering numbers quantify the
complexity of a set). As with other notions touching on convexity, an
important role is played by considerations involving duality. The
central problem in this area is the 1972 ``duality conjecture for
covering numbers" due to Pietsch which has been originally formulated
in the operator-theoretic context, but which in the present notation
can be~stated~as
\begin {cnj} \label{DC}
Do there exist numerical constants $a, b\ge 1$
such that for any dimension $n$
and for any two symmetric convex bodies $K, B$ in $\Rn$ one has
\begin{equation} \label{dc}
      b^{-1}  \log{N(B^{\circ}, a K^{\circ})} \le \log{N(K,B)}
\le b \log{N(B^{\circ}, a^{-1} K^{\circ})} ?
\end{equation}
\end {cnj}
Above and in what follows $A^{\circ} := \{u\in \Rn : \sup_{x\in A}
\langle x, u \rangle \le 1\}$ is the polar body of $A$; ``symmetric" is a
shorthand for  ``symmetric with respect to the origin" and, for
definiteness, all logarithms are to the base~2. In our preferred setting
of a normed space $X$, the proper generality is achieved by considering
$\log{N(K,t \,B)}$ for $t >0$, where $B$ is the unit ball and
$K$ -- a generic (convex, symmetric) subset of $X$. The polars should
then be thought of as subsets of $X^*$, with $B^\circ$
the unit ball of that space. With minimal care, infinite-dimensional
spaces and sets may be likewise considered.  To avoid
stating boundedness/compactness hypotheses,
which are peripheral to the phenomena in question, it is convenient
to allow
$N(\cdot, \cdot)$, $M(\cdot, \cdot)$ etc.  to take the value $+ \infty$.

The quantity $\log{N(K,t \,B)}$ has a clear information-theoretic
interpretation: it is the complexity of $K$, measured in bits,
at the level of resolution $t$ \wrt the metric for which $B$ is the unit
ball. Accordingly, (\ref{dc}) asks whether the complexity of $K$
is controlled by that of the ball in the dual space \wrt
$\|\cdot\|_{K^\circ}$ (the gauge of ${K^\circ}$, i.e., the norm whose unit
ball is ${K^\circ}$), and {\em vice versa},
at {\em every } level of resolution. [The original formulation of the
conjecture involved relating -- in a quantitative way -- compactness of an
operator to that of its adjoint.]
In a very recent paper \cite{AMS} Conjecture \ref{DC} has been verified in
the special yet most important case where
$B$ is an ellipsoid (or, equivalently, when $K$ is a subset of a Hilbert
space);
the reader is referred to that article for a more exhaustive discussion
of historical and mathematical background and for further references.

In the present note we introduce a new notion, which we call
``convex separation" (or ``convexified packing") and prove a duality
theorem related to that concept.  This will lead to a generalization of
the results from \cite{AMS} to the setting requiring only mild geometric
assumptions about the underlying norm. [Both the definition and the generalization are motivated by an earlier paper \cite{BPST}.]
 For example, we now know that
Conjecture \ref{DC} holds -- in the sense indicated in the paragraph
following (\ref{dc}) -- in all $\ell_p$- and $L_p$-spaces (classical or
non-commutative) for
$1<p<\infty$, with constants $a, b$  depending only on $p$ and uniformly
bounded if $p$ stays away from 1 and  $\infty$, and similarly in all
uniformly convex  and all uniformly smooth spaces.

At the same time, and perhaps more importantly, the new
approach ``demystifies" duality results for usual covering/packing
numbers in that it splits the proof into two parts.
One step is a duality theorem for convex separation, which predictably is a
consequence of the  Hahn-Banach theorem. The other step involves geometric
considerations relating convex separation to the usual separation;  while
often delicate and involved, they are always set in a given normed
space and reflect properties of that space without appealing to duality,
and thus are conceptually simpler.

In the next two sections we shall give the definition of convex separation
and prove the corresponding duality theorem. Then, in section 4, we shall
state several estimates for the convexified packing/convex separation
numbers, in particular those that compare them to the usual packing/covering
numbers. In section 5 we state the generalization of the duality result from
\cite{AMS} alluded to above, and give some hints at its proof. We
include details only for the proof of Theorem \ref{tild} (duality for
convex separation) which is the part we consider conceptually new.
Although we see our Theorem \ref{Kconvex} as an essential progress
towards Conjecture \ref{DC}, in this announcement we omit
the proofs; while non-trivial and technically involved, they
require tools which were developed and used in our previous papers.

\bigskip\noindent
{\bf 2. Defining convex separation.} The following notion plays a central
role in this note.  For a set $K$ and a symmetric convex body $B$ we define
   \begin{eqnarray}
      \Nt(K, B)&:=& \sup \{N: \exists \ x_{1}, \ldots, x_{N}\in K \quad
\hbox{ such that } \nonumber \\
      && \qquad \
      (x_{j}+{\rm int}B)\cap {\rm conv}\{x_{i}, i<j\}= \emptyset \},
\label{Mhat}
      \end{eqnarray}
where ``${\rm int}$" stands for the
interior of a set.\footnote {When
defining packing in convex geometry, it is customary to require that only
the {\em interiors} of the translates of
$B$ be disjoint; we follow that convention here even though it is
slightly unsound in the categorical sense.}
We shall refer to any sequence satisfying the condition (\ref{Mhat}) as
$B$-{\em convexly separated}. Leaving out the convex hull operation
``${\rm conv}$" leads to the usual $B$-separated set, which is the same as
$B/2$-packing. Thus we have $\Nt (K, B) \le M(K, B/2)$.
  We emphasize that, as opposed to the usual
notions of packing and covering, the {\em order} of the points is
important here.

The definition (\ref{Mhat}) is very natural from the point of view of
   complexity theory and optimization. A standard device in
constructing geometric algorithms is a ``separation oracle"
(cf. \cite{GLS}): if $T$ is a convex set then, for a given $x$, the oracle
either attests that $x \in T$ or returns a functional efficiently separating
$x$ from $T$. It is arguable that quantities of the type
$\Nt(T,\cdot)$ correctly describe complexities of the set $T$
\wrt many such algorithms.

Since, as pointed out above and in the preceding section, packings, coverings
and separated sets are very closely connected, and the corresponding ``numbers"
are related via two sided estimates involving (at most) small numerical
constants,  in what follows we shall use all these terms interchangeably.

\bigskip\noindent {\bf 3. Duality for $\Nt (K, B)$.}  While it is still an
open problem whether Conjecture \ref{DC} holds in full generality, the
corresponding duality statement for convex separation is fairly
straightforward.
\begin{thm}\label{tild}
      For any pair of convex symmetric bodies $K, B \subset \Rn$ one has
      \[ \Nt (K, B) \le \Nt (B^{\circ}, K^{\circ}/2)^2. \]
      \end{thm}

      \noindent{\em Proof of Theorem \ref{tild}~}
Let $R := \sup \{\|x\|_B : x \in K\}$; i.e., $R$ is the radius of
$K$ \wrt the gauge of $B$.
We will show that
\newline (i)  \, $\Nt (B^{\circ}, K^{\circ}/2) \ge \Nt (K,
B)/\lceil4R\rceil$
\newline (ii) \  $\Nt (B^{\circ}, K^{\circ}/2) \ge \lfloor 4R \rfloor$ + 1

Once the above are proved, Theorem \ref{tild} readily follows.
To show (i), denote $N = \Nt (K, B)$ and let $x_{1},\ldots ,x_{N}$ be a
$B$-convexly separated sequence in  $K$. Then,
by (the elementary version of) the Hahn-Banach theorem, there exist
separating functionals  $y_{1}, \ldots, y_{N} \in B^{\circ}$ such~that
\begin{equation}\label{*}
1 \le i<j \le N \Rightarrow \langle y_{j}, x_{j} - x_{i} \rangle
= \langle y_{j}, x_{j} \rangle -
\langle y_{j}, x_{i} \rangle \ge  1,
\end{equation}
a condition which is in fact equivalent to $(x_j)$ being $B$-convexly
separated. Now  $x_{j} \in K \subset RB$ and
$y_{j}\in B^{\circ}$
imply that $-R \le \langle y_{j}, x_{j} \rangle \le R$, and hence
dividing $[-R,R]$ into $\lceil4R\rceil$ subintervals of length $
\le 1/2$ we may deduce that one of these subintervals contains  $M \ge
N/\lceil4R\rceil$  of the numbers $\langle y_{j}, x_{j} \rangle$.
To simplify the notation, assume that this occurs for $j=1, \ldots , M$,
that is
\begin{equation}\label{**}
1 \le  i,j \le M  \Rightarrow
-1/2 \le \langle y_{i}, x_{i} \rangle -
\langle y_{j}, x_{j} \rangle \le 1/2.
\end{equation}
Combining (\ref{*}) and (\ref{**}) we obtain for $1 \le i<j \le M$
$$
     \langle y_{i}- y_{j} , x_{i} \rangle
= \langle y_{i}, x_{i} \rangle -
\langle y_{j}, x_{j} \rangle + \langle y_{j}, x_{j} \rangle -
\langle y_{j}, x_{i}  \rangle \ge -1/2 +1 =1/2,
$$
which is again a condition of type (\ref{*}) and thus shows that the
sequence $y_{M}, \ldots , y_{1}$ (in this order!) is $K^{\circ}/2$-convexly
separated. Hence
$\Nt( B^{\circ} , K^{\circ}/2) \ge M \ge
N/\lceil4R\rceil$, which is exactly the conclusion of (i).

To show (ii) we note that $R$ is also the radius of $B^\circ$
with respect to the gauge of $K^\circ$. Since we are in a
finite-dimensional space, that radius is attained and so there is a segment
$I :=[-y, y] \subset B^\circ$ with $\|y\|_{K^\circ} =R$.
This implies that $M(I, K^\circ/2) \ge \lfloor 4R \rfloor + 1$. However,
in dimension one separated and convexly separated sets coincide;
this allows to conclude that
$\Nt(B^\circ, K^\circ/2) \ge \Nt(I, K^\circ/2) \ge \lfloor 4R \rfloor + 1$,
as required.

\bigskip\noindent {\bf 4. Separation vs. convex separation.}  It appears at
the first sight that the notion of convex separation is much more
restrictive than that of usual separation and, consequently, that --
except for very special cases such as that of one-dimensional sets
mentioned above --  $\Nt(\cdot, \cdot)$ should be significantly smaller
than
$M(\cdot, \cdot)$ or $N(\cdot, \cdot)$. However, we do not have examples
when that happens. On the other hand, for several interesting classes of
sets we do have equivalence for not-so-trivial reasons.
Here we state two such results.

\begin{thm}\label{ellipsoids}
      There exist numerical constants $C, c >0$ such that, for any
$n \in \mathbb{N}$ and any pair of ellipsoids $\mathcal{E}, \mathcal{B}
\subset
\Rn$ one has
      \[\log M (\mathcal{E}, \mathcal{B}) \le
C\log \Nt (\mathcal{E}, c \,\mathcal{B}).
\]
      \end{thm}

While Theorem \ref{ellipsoids} deals with purely Euclidean setting, the next
result holds under rather mild geometric assumptions about
the underlying norm. It requires $K$-convexity, a property which goes back to
\cite{Mau-Pis} and is well known to specialists;
we refer to \cite{Pi} for background and properties. While many interesting
descriptions of that class are possible, here we just mention that
$K$-convexity
is equivalent to the absence of large subspaces well-isomorphic to
finite-dimensional
$\ell_1$-spaces and that it can be quantified, i.e., there is a
parameter called
the $K$-convexity constant and denoted $K(X)$, which can be defined 
both for 
finite and infinite 
dimensional spaces, and which has good permanence properties
with respect
to standard functors of functional analysis. For example, as hinted in the
introduction, all $L_p$-spaces for $1<p<\infty$ are
$K$-convex (with constants  depending only on $p$), and similarly all
uniformly convex and all uniformly smooth spaces.

\begin{thm}\label{boundedtype} Let $X$ be a normed space which is
$K$-convex with $K(X) \le \kappa$ and let $B$ be its unit ball.  Then for any
  bounded  set
$T\subset X$ one has
$$
\log M(T, B) \le \beta \log \Nt  (T, B/2),
$$
where  $\beta$ depends only on $\kappa$ and the diameter of $T$. Similarly, if
$r>0$ and if $U$ is a symmetric convex subset of $X$ with $U \supset rB$,
then
$$
   \log M(B, U) \le \beta' \log \Nt  (B, U/2)
$$
with  $\beta'$ depending only on $\kappa$ and $r$.
      \end{thm}

The proof of Theorem \ref{ellipsoids} is non-trivial but elementary.  The
proof of Theorem \ref{boundedtype} is based on the so called Maurey's lemma
(see \cite{Pi0}) and the ideas from \cite{BPST}.  The details of both
arguments will be presented elsewhere.

\bigskip\noindent {\bf 5. Duality of covering and packing numbers in
$K$-convex spaces.}  If $B$ is the unit ball in a $K$-convex space $X$,
then, combining Theorems   \ref{tild} and \ref{boundedtype}, we obtain for any
bounded symmetric convex set $T \subset X$
$$
\log M(T, B) \le \beta \log \Nt  (T, B/2) \le 2 \beta\log \Nt
(B^{\circ}, T^{\circ}/4) \le 2 \beta\log M(B^{\circ}, T^{\circ}/8),
$$
where $\beta$ depends only on the diameter of $T$ (and on $K(X)$) and,
similarly,
$\log M(B^{\circ}, T^{\circ}) \le 2 \beta' \log M(T, B/8)$.  To show the
latter, we apply the second part of Theorem \ref{boundedtype} to
$X^*$, the dual
of $X$, and to $U=T^{\circ} \supset (2/{\rm diam} \, T) B^\circ$ (see
the comments following (\ref{dc})), and use
the known fact that $K(X^*)=K(X)$.

On the other hand, an iteration scheme developed in
\cite{AMS} can be employed to show that if, for some normed space $X$, duality
in the sense of the preceding paragraph holds -- with some constants
$\beta, \beta'$ -- for, say,  all $T \subset 4B$, then it also holds for {\em
all } sets $T \subset X$ with constants depending only on $\beta, \beta'$.
We thus have
\begin{thm}\label{Kconvex} Let $X$ be a normed space which is
$K$-convex with $K(X) \le \kappa$ and let $B$ be its unit ball.  Then for any
symmetric convex set $T\subset X$ and any $\epsilon >0$ one has
$$
b^{-1} \log M  (B^\circ, a \epsilon T^\circ) \le
\log M(T, \epsilon B) \le b \log M  (B^\circ, a^{-1} \epsilon T^\circ),
$$
where  $a, b \ge 1$ depend only on $\kappa$.
      \end{thm}
Theorem \ref{Kconvex} is related to \cite{AMS} in a similar way as
the results of \cite{BPST} were related to \cite{T}:  in both cases
a statement concerning duality of covering numbers is generalized from
the Hilbertian setting to that of a $K$-convex space. A more detailed
exposition
of its proof will be presented elsewhere.

\bigskip\noindent {\bf 6. Final remarks.}  While the approach via convexly
separated sets does not yet settle Conjecture
\ref{DC} in full generality,
it includes, in particular,
all cases for which the Conjecture has been previously
verified.  One such special case, which does not follow directly from
the results included in the preceding sections,  was settled in \cite{KoM}
and subsequently generalized in
\cite{Pi1}:   {\em If, for some } $\gamma > 0$,
   $\log{N(K,B)} \ge \gamma n$, {\em then } $\log{N(B^{\circ}, K^{\circ}/a)} \ge
b^{-1} \log{N(K,B)}$ {\em with } $a, b \ge 1$ {\em depending only on
} $\gamma$.
In the same direction, we have
\begin{thm}\label{expo}
Let $\gamma > 0$ and let $K, B \subset \Rn$ be symmetric convex bodies with
$\log M(K, B) \ge \gamma n$. Then
\begin{equation} \label{expoeq}
\log \Nt (K, B/\alpha) \ge \beta^{-1} \log M(K, B),
\end{equation}
where $\alpha \ge 1$ depends only on $\gamma$ and $\beta \ge 1$ is a universal
constant.
\end{thm}
The result from \cite{KoM}, \cite{Pi1} mentioned above is then an easy
corollary: combine (\ref{expoeq}) with Theorem \ref{tild}  to obtain
$\log M(K, B) \le 2 \beta \log \Nt(B^{\circ}, (2 \alpha)^{-1} K^{\circ})$,
and the later expression is ``trivially"
$\le 2 \beta \log M(B^{\circ}, (4 \alpha)^{-1} K^{\circ})$.  In turn, Theorem
\ref{expo} can be derived by applying Theorem \ref{ellipsoids} to the
so-called $M$-ellipsoids of $K$ and $B$.  This is admittedly not the simplest
argument, but it subscribes to
our philosophy of minimizing the role of duality considerations and sheds
additional light on the relationship between $M(\cdot, \cdot)$ and $\Nt (\cdot,
\cdot)$.

We conclude the section and the note by pointing out that
a slightly different -- but equally natural from the algorithmic point of
view -- definition of convex separation is possible. Namely, we may
consider sets of points in $K$ verifying (\ref{Mhat}) with a modified condition
$(x_{j}+{\rm int}B)\cap {\rm conv}\{x_{i}, i \neq j\}= \emptyset$.
For the so modified convex separation number -- let us call it
$\tilde{M}(K, B)$ --
the statement (i) from the proof of Theorem \ref{tild} remains true
(by an almost identical argument), and so for ``bounded" sets we do
have duality.
However, the statement (ii) is false already for the simple case of
a segment.

\small
\bigskip\noindent {\em Acknowledgement.}
{
This research has been supported in part by grants from
the U.S.-Israel Binational Science Foundation for the first three named
authors and from the National Science Foundation (U.S.A.)
for the third-named author.
The fourth-named author holds the Canada
Research Chair in Geometric  Analysis.
The first and second-named authors thank respectively
University Paris VI and IHES  for their hospitality
during the time when this work was in progress.

\bigskip
{\scriptsize
\noindent{S. Artstein, V. D. Milman

\noindent School of Mathematical Sciences,
Tel Aviv University,
Tel Aviv 69978, Israel

\noindent E-mail: artst@post.tau.ac.il, milman@post.tau.ac.il

\vskip 8pt

\noindent S. J. Szarek

\noindent Department of Mathematics,
Case Western Reserve University,Cleveland, OH 44106-7058, U.S.A.

{\em and }

\noindent Equipe d'Analyse Fonctionnelle, B.C. 186,
Universit\'{e} Paris VI,
4 Place Jussieu,  75252  Paris, France

\noindent E-mail: szarek@cwru.edu

\vskip 8pt

\noindent N. Tomczak-Jaegermann

\noindent Department of Mathematical and Statistical Sciences, University of
Alberta,

\noindent Edmonton, Alberta, Canada T6G 2G1

\noindent E-mail:   nicole@ellpspace.math.ualberta.ca

}}

\end {document}